\documentclass{article}

\usepackage{amsmath}
\usepackage{amssymb}

\newtheorem{thm}{Theorem}[section]

\newtheorem{cor}[thm]{Corollary}

\newtheorem{lem}[thm]{Lemma}
\newtheorem{rem}[thm]{Remark}
\newtheorem{ex}[thm]{Example}

\newcommand{\gp}{\mathfrak{p}}
\newcommand{\gq}{\mathfrak{q}}
\newcommand{\gm}{\mathfrak{m}}
\newcommand{\gn}{\mathfrak{n}}
\newcommand{\ga}{\mathfrak{a}}
\newcommand{\gb}{\mathfrak{b}}
\newcommand{\gc}{\mathfrak{c}}

\newcommand{\pp}{\mathcal{P}}
\newcommand{\zz}{\mathbb{Z}}

\newcommand{\ic}[1]{\overline{#1}}

\newcommand{\mult}[2]{{\mathrm e}_{#1}(#2)}

\newcommand{\lc}[3]{{\mathrm H}^{#1}_{#2}(#3)}

\newcommand{\jmult}[3]{{\mathrm j}_{#1}(#2, #3)}
\newcommand{\jm}[1]{{\mathrm j}\,(#1)}
\newcommand{\length}[2]{{\rm length}_{\,#1}\,#2}
\newcommand{\ass}[2]{{\rm Ass}_{#1}\,(#2)}
\newcommand{\assh}[2]{{\rm Assh}_{#1}\,(#2)}
\newcommand{\Min}[2]{{\rm Min}_{#1}\,(#2)}
\newcommand{\height}[2]{{\rm ht}_{#1}\,#2}

\newcommand{\ann}[2]{{\rm ann}_{#1}\,#2}
\newcommand{\kdim}[2]{{\rm dim}_{#1}\,#2}

\newcommand{\ra}{\longrightarrow}
\newcommand{\homom}[1]{\stackrel{#1}{\longrightarrow}}

\def\trmatrix#1{\sideset{^t\!\!}{}{\mathop{\,#1}}} 

\begin{document}

\title{\Large Computing j-multiplicity}

\author{$\mbox{Koji Nishida}^\ast$
\and 
$\mbox{Bernd Ulrich}^\dagger$}
\date{
\begin{flushleft}
\small
\hspace{10ex} ${}^\ast$\,Department of Mathematics and Informatics, Chiba University, \\
  \hspace*{15ex} 1-33 Yayoi-Cho, Inage-Ku, Chiba-Shi, 263-8522, Japan; \\
  \hspace*{15ex} e-mail: nishida@math.s.chiba-u.ac.jp \\
\hspace{10ex} ${}^\dagger$\,Department of Mathematics, Purdue University, \\
  \hspace*{15ex} West Lafayette, IN 47907-1395, USA; \\
  \hspace*{15ex} e-mail: ulrich@math.purdue.edu
\end{flushleft}
}
\maketitle

\begin{abstract}
The j-multiplicity is an invariant that can be defined for
any ideal in a Noetherian local ring $(R, \gm)$.
It is equal to the Hilbert-Samuel multiplicity
if the ideal is $\gm$-primary.
In this paper we explore the computability of the j-multiplicity,
giving another proof for the length formula and
the additive formula.
\end{abstract}

\section{Introduction}
Let $(R, \gm)$ be a Noetherian local
ring of of dimension $d > 0$ and let $I$ be a proper
ideal of $R$.
We set
\[
\jm{I} = \lim_{n \rightarrow \infty}\,
\frac{(d - 1)!}{n^{d-1}}\hspace{1ex}
\length{R}{\lc{0}{\gm}{I^n/I^{n+1}}}
\]
and call it the {\it j-multiplicity} of $I$,
where $\lc{0}{\gm}{\hspace{1ex}\cdot\hspace{1ex}}$
denotes the $0$-th local cohomology functor.
It is obvious that $\jm{I}$ coincides with
the usual multiplicity $\mult{I}{R}$ if $I$
is an $\gm$-primary ideal.
Moreover, the j-multiplicity enjoys a lot of
properties similar to those of the usual multiplicity
(see \cite{am}, \cite{fm1}, \cite{fm2}, \cite{fov}).

In this paper we show some computation of the j-multiplicity
of certain ideals using the length formula which was first
proved by Achilles and Manaresi in \cite[(3.8)]{am}.
We will give another proof for the formula
taking a different approach.
The length formula can be stated as follows:
If $A/\gm$ is an infinite field,
then choosing sufficiently general elements
$a_1, \dots\,, a_{d-1}, a_d$ of $I$, we have
\[
\jm{I} =
\length{R}{R\, /\, {((a_1, \dots\,, a_{d-1}) :_R I^n)+a_dR}}
\]
for $n \gg 0$.
Let $\ga = (a_1, \dots\,, a_{d-1})$ and
$\gb = \ga :_R I^n$ for $n \gg 0$.
Then the equality above means
$\jm{I} = \mult{a_dR}{R / \gb} =
\mult{I}{R / \gb}$
since $R / \gb$ is the zero ring or
a Cohen-Macaulay ring of dimension one.
If $R$ itself is Cohen-Macaulay and
$I$ satisfies a certain Artin-Nagata property
(cf. \cite{u}),
then $I \cap \gb = \ga$, and so we have an exact sequence
$0 \rightarrow I / \ga \rightarrow
R / \gb \rightarrow R / I + \gb \rightarrow 0$,
which implies
$\mult{a_dR}{R / \gb} = \mult{a_dR}{I / \ga} =
\length{R}{I / \ga + a_dI}$\,.
Thus we get
\[
\jm{I} = \length{R}{I / \ga + a_dI}\,,
\]
which is the length formula of Flenner and Manaresi
proved in \cite{fm2}.

Section 2 of this paper is devoted to recalling
some known facts about j-multiplicity in the graded case.
Our proofs are slightly simpler than the original ones.
In Section 3, we prove the length formula
and give some corollaries.
Applying our argument, we will give yet another proof
for the additivity of the j-multiplicity in
\ref{3.8}.
Moreover, we consider the associativity formula
of j-multiplicity in \ref{3.10}.
In the last section, we investigate some concrete
examples in order to illustrate how to use
the formulas given in section 3.
\section{Graded case}
Let $G = \oplus_{n \geq 0}\,G_n$ be a positively
graded Noetherian ring such that
$(G_0, \gn)$ is local and $G = G_0[\,G_1\,]$.
We assume that $G_0 / \gn$ is infinite.
Let $L = \oplus_{n \geq 0}\,L_n$
be a finitely generated graded $G$-module.
We set $W = \lc{0}{\gn G}{L}$ and take an integer $k \gg 0$.
Then, as $W \cap \gn^kL = 0$,
$W$ can be imbedded into $L / \gn^k L$,
and so $\kdim{G}{W} \leq \ell(G, L)$,
where $\ell(G, L)$ denotes the Krull dimension of
$L / \gn L$ as a $G$-module.
Furthermore, as $W$ is a finitely generated graded
$G / \gn^k G$-module and as $G / \gn^k G$
is a standard algebra over an Artinian local ring,
there exists a polynomial $P(n)$ in $n$
such that $\length{G_0}{W_n} = P(n)$ for $n \gg 0$.
Now we take a positive integer $d$ such that
$\kdim{G}{L} \leq d$.
Then the degree of $P(n)$ is at most $d-1$
and the coefficient of $n^{d-1}$ is
$\alpha / (d-1)!$ for some integer $\alpha \geq 0$.
We denote this integer $\alpha$ by $\jmult{d}{G}{L}$.
Hence we have
\[
\jmult{d}{G}{L} =
\lim_{n \rightarrow \infty}
\frac{(d-1)!}{n^{d-1}}\hspace{1ex} \length{G_0}{W_n} \geq 0\,.
\]
It is obvious that
$\jmult{d}{G}{L} \neq 0$ if and only if
$\kdim{G}{W} = d$.
Furthermore, we have the following.
\begin{lem}\label{2.1}
{\rm (\cite[6.1.6]{fov})}\hspace{1ex}
$\jmult{d}{G}{L} \neq 0$ if and only if
$\ell(G, L) = d$.
\end{lem}
{\it Proof.}\hspace{0.5ex}
Suppose $\ell(G, L) = d$.
Then there exists $P \in {\rm Supp}_{\,G}\,L / \gn L$ such that
$\dim\, G / P = d$.
Because $P \in {\rm Supp}_{\,G}\,L$ and
$\kdim{G}{L} \leq d$, we have $P \in \ass{G}{L}$.
Hence $P = \ann{G}{x}$ for some $x \in L$.
We notice that $\gn \cdot x = 0$ as $\gn G \subseteq P$.
This means $x \in W$,
and so $P \in \ass{G}{W}$.
Thus we get $\kdim{G}{W} = d$.
The converse is obvious.

\vspace{1em}
Let $f \in G$.
We say that $f$ is $G_+$-filter regular for $L$
if $f \not\in P$ for any $P \in \ass{G}{L}$ with
$G_+ \not\subseteq P$.
In the case where $f \in G_1$,
it is well known that $f$ is
$G_+$-filter regular for $L$ if and only if
the multiplication map
$L_{n-1} \homom{f} L_n$
is injective for $n \gg 0$.
\begin{lem}\label{2.2}
Let $f \in G_1$ be $G_+$-filter regular for $L$.
Then we have the following assertion.
\begin{itemize}
\item[{\rm (i)}]
$f$ is $G_+$-filter regular for $L/W$.
\item[{\rm (ii)}]
$fL_{n-1} \cap W_n = fW_{n-1}$ for $n \gg 0$.
\item[{\rm (iii)}]
If $\kdim{G_0}{L_n} < \kdim{G}{L}$ for all $n \geq 0$,
then $f$ is an ssop for $L$.
\end{itemize}
\end{lem}
{\it Proof.}\hspace{0.5ex}
Because $W = 0 :_{L} \gn^k G$ for some $k > 0$,
we have $\ass{G}{L / W} \subseteq \ass{G}{L}$,
which implies (i).
Then we get (ii) since the multiplication map
$L_{n-1} / W_{n-1} \homom{f}
L_n / W_n$ is injective for $n \gg 0$.
In order to prove (iii),
take $r \gg 0$ so that $L$ can be
generated by the elements in $L_0, L_1, \dots\,, L_r$
and put $L' = \oplus_{n > r}\,L_n$.
Then
\[
G_0 \cap (L' :_G L) \subseteq
\cap_{n = 0}^r\,\ann{G_0}{L_n} \subseteq \ann{G}{L}\,,
\]
and so $G_0 \cap (L' :_G L) \subseteq Q$
for any $Q \in {\rm Supp}_{\,G}\,L$.
Hence, if $Q \in {\rm Supp}_{\,G}\,L$ contais $G_+$,
we have $L' :_G L \subseteq Q$,
and so $Q \in {\rm Supp}_{\,G}\,L / L'$.
This implies that if
$\assh{G}{L} \cap {\mathrm V}(G_+) \neq \phi$,
then $\dim_{\,G}\,L / L' = \dim_{\,G}\,L$,
which implies $\dim_{\,G_0}\,L_n =
\dim_{\,G}\,L$ for some $0 \leq n \leq r$.
Thus we get (iii).
\begin{lem}\label{2.3}
Let $d \geq 2$ and $\kdim{G_0}{L_n} < \kdim{G}{L}$
for any $n \geq 0$.
Then, choosing $f \in G_1$ generally, we have
$\kdim{G}{L / fL} \leq d - 1$ and
$\jmult{d}{G}{L} = \jmult{d-1}{G}{L / fL}$.
\end{lem}
{\it Proof.}\hspace{0.5ex}
Let us choose $f \in G_1$ so that $f$ is $G_+$-filter
regular for $L$.
Then $\kdim{G}{L / fL} \leq d - 1$ by (iii) of \ref{2.2}.
Furthermore, as $\gn G$ contains a non-zero-divisor
of $L / W$, we have $\ell(G, L / W) \leq d - 1$,
and hence we can choose $f \in G_1$ so that
$\ell(G, L / fL + W) \leq d - 2$.
This means $\jmult{d-1}{G}{L / fL + W} = 0$ by \ref{2.1}.
Now we take $n \gg 0$.
Then $fL_{n-1} \cap W_n = fW_{n-1}$ by (ii) of \ref{2.2},
and so we have an exact sequence
\[
0 \ra W_n / fW_{n-1} \ra L_n / fL_{n-1}
\ra L_n / {fL_{n-1}+W_n} \ra 0
\]
Because $W_n / fW_{n-1}$ is a $G_0$-module of finite length,
we get the exact sequence
\[
0 \ra W_n / fW_{n-1} \ra \lc{0}{\gn}{L_n / fL_{n-1}}
\ra \lc{0}{\gn}{L_n / fL_{n-1}+W_n} \ra 0\,.
\]
Therefore $\jmult{d-1}{G}{L / fL} =
\jmult{d-1}{G}{W / fW}$.
On the other hand, as
\[
0 \ra W_{n-1} \homom{f} W_n \ra W_n / fW_{n-1} \ra 0
\]
is exact for $n \gg 0$, we have
\begin{eqnarray*}
\jmult{d-1}{G}{W / fW} & = &
 \lim_{n \rightarrow \infty}
   \frac{(d-2)!}{n^{d-2}}\hspace{1ex}
   \length{G_0}{W_n/fW_{n-1}}  \\
 & = &
 \lim_{n \rightarrow \infty}
   \frac{(d-2)!}{n^{d-2}}\hspace{1ex}
   (\length{G_0}{W_n} - \length{G_0}{W_{n-1}}) \\
 & = &
 \jmult{d}{G}{L}\,,
\end{eqnarray*}
and the proof is complete.
\section{Local case}
Let $I$ be a proper ideal of a Noetherian local ring
$(R, \gm)$ such that $\kdim{}{R} > 0$ and
$R / \gm$ is an infinite field.
Let $M$ be a finitely generated $R$-module and
$d$ be a positive integer such that $\kdim{R}{M} \leq d$.
Then ${\rm gr}_{I}\,(M) = \oplus_{n \geq 0}\,I^nM / I^{n+1}M$
is a finitely generated graded ${\rm gr}_{I}\,(R)$-module
whose Krull dimension is at most $d$.
We set
\[
\jmult{d}{I}{M} := \jmult{d}{{\rm gr}_{I}\,(R)}{{\rm gr}_{I}\,(M)}
\]
and call it the {\it j-multiplicity} of $I$
with respect to $M$. Hence we have
\[
\jmult{d}{I}{M} = \lim_{n \rightarrow \infty}\,
\frac{(d - 1)!}{n^{d-1}}\hspace{1ex}
\length{A}{\lc{0}{\gm}{I^nM/I^{n+1}M}}\,,
\]
and so, if $\kdim{R}{M} = d$ and $M / IM$
has finite length, $\jmult{d}{I}{M}$ coincides with
the usual multiplicity $\mult{I}{M}$.
We denote
$\jmult{\dim R}{I}{R}$ by $\jm{I}$.
Let $\ell(I, M) := \ell({\rm gr}_{I}\,(R), {\rm gr}_{I}\,(M))$
and $\ell(I) = \ell(I, R)$.
\begin{lem}\label{3.1}
$\jmult{d}{I}{M} \neq 0$ if and only if
$\ell(I, M) = d$.
\end{lem}
{\it Proof.}\hspace{0.5ex}
This follows immediately from \ref{2.1}.
\begin{lem}\label{3.2}
Let $d \geq 2$ and $\kdim{R}{M/IM} < \kdim{R}{M}$.
Then, for a general element $a \in I$,
we have $\kdim{R}{M/aM} \leq d - 1$
and $\jmult{d}{I}{M} = \jmult{d-1}{I}{M/aM}$.
\end{lem}
{\it Proof.}\hspace{0.5ex}
We apply \ref{2.3} setting $G = {\rm gr}_{I}\,(R)$ and
$L = {\rm gr}_{I}\,(M)$.
For any $n \geq 0$, $L_n$ is a homomorphic image
of the direct sum of finite number of copies of $M/IM$,
and hence we have
\[
\kdim{G_0}{L_n} \leq \kdim{R}{M/IM} < \kdim{R}{M} = \kdim{G}{L}\,.
\]
Now we choose a general element $a \in I$.
Let $f$ be the initial form of $a$ in $G$.
Then $\kdim{G}{L/fL} \leq d - 1$ and
$\jmult{d}{I}{M} = \jmult{d-1}{G}{L/fL}$ by \ref{2.3}.
Because ${\rm gr}_{I}\,(M/aM)$ is a homomorphic image of
$L/fL$, we have $\kdim{R}{M/aM} \leq d - 1$.
Moreover, we may assume that $a$ is a superficial element of $I$
with respect to $M$,
and so $[L/fL]_n \cong [{\rm gr}_{I}\,(M/aM)]_n$ for $n \gg 0$.
Then we have $\jmult{d-1}{G}{L/fL} = \jmult{d-1}{I}{M/aM}$
by the definition of j-multiplicity and the proof is complete.

\vspace{1em}
For an $R$-submodule $K$ of $M$, we set $K :_M I^\infty =
\bigcup_{n \geq 0} (K :_M I^n)$.
\begin{lem}\label{3.3}
Let $\ic{M} = M / {0 :_M I^\infty}$.
Then $\jmult{d}{I}{M} = \jmult{d}{I}{\ic{M}}$.
Furthermore we have $\kdim{R}{\ic{M}/I\ic{M}} < \kdim{R}{\ic{M}}$
if $\ic{M} \neq 0$.
\end{lem}
{\it Proof.}\hspace{0.5ex}
Let $n \gg 0$.
Then $I^nM \cap (0 :_M I^\infty) = 0$,
which implies $I^n\ic{M}/I^{n+1}\ic{M} \cong I^nM/I^{n+1}M$,
and hence the first assertion follows.
On the other hand, as $I$ contains a non-zero-divisor of $\ic{M}$,
we get the second assertion.
\begin{lem}\label{3.4}
Let $a$ be a general element of $I$.
We put $\ic{M} = M/{(0 :_M I^\infty)}$.
Then we have $\kdim{R}{\ic{M}/a\ic{M}} \leq d - 1$ and
$\jmult{d}{I}{M} = \jmult{d-1}{I}{\ic{M}/a\ic{M}}$.
\end{lem}
{\it Proof.}\hspace{0.5ex}
If $\ic{M} \neq 0$,
the assertion immediately follows from \ref{3.2} and \ref{3.3}.
If $\ic{M} = 0$, then $\jmult{d}{I}{M} = 0$
since $I^nM = 0$ for $n \gg 0$,
and so the assertion is still true.
\begin{lem}\label{3.5}
Let $1 \leq i < d$.
Then, choosing general elements
$a_1, \dots\,, a_i$ of $I$,
we have $\kdim{R}{M/((a_1, \dots\,, a_i)M :_M I^\infty)}
\leq d - i$ and
\[
\jmult{d}{I}{M} =
\jmult{d-i}{I}{M/{(a_1, \dots\,, a_i)M :_M I^\infty}}\,.
\]
\end{lem}
{\it Proof.}\hspace{0.5ex}
We use induction on $i$.

Let us consider the case where $i = 1$.
We set $\ic{M} = M/(0 :_M I^\infty)$ and
$M' = \ic{M}/a_1\ic{M}$.
Then by \ref{3.3} and \ref{3.4} we have $\kdim{R}{M'} \leq d - 1$
and
$\jmult{d}{I}{M} = \jmult{d-1}{I}{M'}
= \jmult{d-1}{I}{M'/(0 :_{M'} I^\infty)}\,$.
Because the kernel of the composition of the canonical surjections
$M \rightarrow \ic{M} \rightarrow M' \rightarrow
M'/(0 :_{M'} I^\infty)$ is $a_1M :_M I^\infty$,
we get $\kdim{R}{M/(a_1M :_M I^\infty)} \leq d - 1$
and
$\jmult{d}{I}{M} = \jmult{d-1}{I}{M/(a_1M :_M I^\infty)}\,$.

Now we suppose $1 < i < d$ and set
$N = M/((a_1, \dots\,, a_{i-1})M :_M I^\infty)$.
Then $\kdim{R}{N} \leq d - i +1$ and
$\jmult{d}{I}{M} = \jmult{d-i+1}{I}{N}$ by the inductive hypothesis.
Furthermore, by the assertion in the case where $i = 1$,
we have $\kdim{R}{N/(a_iN :_N I^\infty)} \leq d - i$
and $\jmult{d-i+1}{I}{N} =
\jmult{d-i}{I}{N/(a_iN :_N I^\infty)}$.
on the other hand, the kernel of the composition
of the canonical surjections
$M \rightarrow N \rightarrow N/(a_iN :_N I^\infty)$
is $(a_1, \dots\,, a_i)M :_M I^\infty$.
Therefore we get the required assertion.
\begin{thm}\label{3.6}
Choosing general elements
$a_1, \dots\,, a_{d-1}, a_d$ of $I$, we have
\begin{eqnarray*}
\jmult{d}{I}{M} & = & 
 \mult{I}{M\,/\,((a_1, \dots\,, a_{d-1})M :_M I^\infty)} \\
      & = & 
 \length{R}{M\,/\,{((a_1, \dots\,, a_{d-1})M :_M I^\infty) + a_dM}}\,. 
\end{eqnarray*}
\end{thm}
{\it Proof.}\hspace{0.5ex}
We put $N = M/{((a_1, \dots\,, a_{d-1})M :_M I^\infty)}$.
Then $\kdim{R}{N} \leq 1$ and $\jmult{d}{I}{M} =
\jmult{1}{I}{N}$ by \ref{3.5}.
Therefore it is enough to show that
$\jmult{1}{I}{N} = \length{R}{N/{a_dN}}$.
This is obvious if $N = 0$.
So, let us consider the case where $N \neq 0$.
In this case, we may assume that $a_d$ is an $N$-regular
element and $a_dR$ is a reduction of $I$ with respect to $N$.
Moreover, $N/IN$ has finite length.
Therefore we have $\jmult{1}{I}{N} =
\mult{I}{N} = \mult{a_dR}{N} = \length{R}{N/{a_d}N}$
and the proof is complete.
\begin{rem}\label{3.11}
Suppose that $I$ is generated by $d$ elements
and $b \in I \setminus \gm I$.
Then we can choose general elements
$a_1, \dots\,, a_{d-1}$ of $I$
so that they generate $I$ together with $b$.
In this case we have
\[
\jmult{d}{I}{M} =
\length{R}{M\,/\,((a_1, \dots\,, a_{d-1})M :_M I^\infty}) + bM\,.
\]
\end{rem}

Let $S$ be a Noeterian ring and $s$ an integer.
We say that an ideal $J$ of $S$
satisfies $G_s$
if $\mu_{S_P}(JS_P) \leq \height{S}{P}$
for any $P \in {\rm V}(J)$ with $\height{S}{P} < s$.
An ideal $K$ of $S$ is called an
$s$-residual intersection of $J$
if $K = \ga :_S J$ for some $s$-generated ideal $\ga \subseteq J$
with $\height{S}{(\ga :_S J)} \geq s$.
If $K$ is an $s$-residual intersection of $J$ with
$\height{S}{(J + K)} \geq s + 1$,
then $K$ is called a geometric $s$-residual intersection of $J$.
\begin{cor}\label{3a}
Let $S = \oplus_{n \geq 0} S_n$ be a
$d$-dimensional standard algebra such that $S_0 = k$
is an infinite field.
We put $\gn = S_+$ and $R = S_\gn$.
Let $0 < r \in \zz$ and $J$ be an ideal of $S$
generated by homogeneous elements of degree $r$.
We put $I = JR$.
Suppose that $\ell(I) = d$ and
$J$ satisfies $G_d$.
Then, if $f_1, \dots\,, f_{d-1} \in J$ are
general homogeneous elements of degree $r$, we have
\[
\jm{I} = r \cdot
\mult{}{S\,/\,((f_1, \dots\,, f_{d-1}) :_S J)}\,,
\]
where $\mult{}{\ast}$ denotes the multiplicity
with respect to $\gn$.
\end{cor}
{\it Proof.}\hspace{0.5ex}
We put $\ga = (f_1, \dots\,, f_{d-1})$ and
$T = S\,/\,(\ga :_S J^\infty)$.
As $\ell(I) = d$, we have $T \neq 0$,
and so $T$ is a Cohen-Macaulay $S$-module
with $\kdim{S}{T} = 1$.
We notice that $\mult{}{T} = \length{k}{T_i}$ for $i \gg 0$.
On the other hand, by \ref{3.6} it follows that
$T\,/\,JT$ has finite length and
$\jm{I} = \mult{J}{T}$.
Now we take a homogeneous element $f_d \in J$
of degree $r$ so that $f_dS$ is a reduction of $J$
with respect to $T$.
Then $\jm{I} = \length{S}{T\,/\,f_dT}$.
From the exact sequence
$0 \ra T(-r) \stackrel{f_d}{\ra}
T \ra T\,/\,f_dT \ra 0$
of graded $S$-modules, we get
$\length{k}{[T\,/\,f_dT]_i} =
\length{k}{T_i} - \length{k}{T_{i-r}}$
for any $i \in \zz$.
Therefore, taking $n \gg 0$, we have
\[
\length{S}{T\,/\,f_dT} = \sum_{i = n-r+1}^n T_i =
r \cdot \mult{}{T}\,.
\]
Thus we get $\jm{I} = r \cdot \mult{}{T}$,
and hence it is enough to show that
$\mult{}{T} = \mult{}{T'}$,
where $T' = S\,/\,\ga :_S J$.
Because $J$ satisfies $G_d$,
$J + (\ga :_S J)$ is $\gn$-primary.
Hence $(\ga :_S J^\infty)\,/\,(\ga :_S J)$
has finite length as it is annihilated by
$J^\nu + (\ga :_S J)$ for $\nu \gg 0$.
Therefore we get $\mult{}{T} = \mult{}{T'}$
from the exact sequence
$0 \ra (\ga :_S J^\infty)\,/\,(\ga :_S J) \ra T' \ra T \ra 0$,
and the proof is complete.

\vspace{1em}

We say that an ideal $J$ in a Noetherian ring $S$ is
({\it weakly}) {\it $s$-residually $S_2$}
if for every $i$ with $\height{S}{J} \leq i \leq s$
and every (geometric) $i$-residual intersection $K$ of $J$,
$R\,/\,K$ is $S_2$.
\begin{cor}\label{3b}
Let $R$ be Cohen-Macaulay and $d = \kdim{}{R}$.
Suppose that $I$ satisfies $G_d$ and $I$ is
weakly $(d-2)$-residually $S_2$.
Then, choosing general elements
$a_1, \dots\,, a_{d-1}, a_d$ of $I$, we have
\[
\jm{I} = \length{R}{R\,/\,((a_1, \dots\,, a_{d-1}) :_R I) + a_dR}\,.
\]
\end{cor}
{\it Proof.}\hspace{0.5ex}
We put $\ga = (a_1, \dots\,, a_{d-1})$ and $N = R\,/\,(\ga :_R I)$.
As $I$ satisfies $G_d$, $\ga :_R I$ is a geometric
$(d - 1)$-residual intersection of $I$ by \cite[1.4]{u},
and so $N$ is the zero module or
a Cohen-Macaulay $R$-module with $\kdim{R}{N} = 1$
(cf. \cite[3.4]{ceu}).
In particular, $\ga R_\gp = IR_\gp$
for any $\gp \in {\rm V}(I) \setminus \{ \gm \}$
and $\gm \not\in \ass{R}{N}$.
Hence $\ga :_R I = \ga :_R I^\infty$, and so
$\jm{I} = \mult{I}{N} =
\length{R}{N\,/\,a_dN}$ by \ref{3.6},
which means the required equality.

\vspace{1em}

Next we give another proof for the additivity of j-multiplicity
which was first proved by Flenner, O'Carroll and Vogel in 
\cite[6.1.7]{fov}.
For that purpose, we prepare the following result.

\begin{lem}\label{3.7}
Let $d = 1$.
We put $\pp = \{ \gp \in {\rm Spec}\,R \mid
\mbox{$\kdim{}{R/\gp} = 1$ and $I \not\subseteq \gp$}\}$.
Then
\[
\jmult{1}{I}{M} = \sum_{\gp \in \pp}\,
\length{R_\gp}{M_\gp}\cdot \mult{I}{R/\gp}\,.
\]
\end{lem}
{\it Proof.}\hspace{0.5ex}
As $\jmult{1}{I}{M} = \jmult{1}{I}{M/(0 :_M I^\infty)}$
by \ref{3.3} and as $M_\gp = (M/(0 :_M I^\infty))_\gp$
for any $\gp \in \pp$,
we may assume ${0 :_M I^\infty} = 0$ and $M \neq 0$.
Then $I$ contains an $M$-regular element.
Hence $\kdim{R}{M} = 1$ and $M/{IM}$ has finite length.
Therefore we have
\[
\jmult{1}{I}{M} = \mult{I}{M} = \sum_{\gp \in \assh{R}{M}}\,
\length{R_\gp}{M_\gp}\cdot \mult{I}{R/\gp}\,.
\]
As $I$ contains an $M$-regular element,
$I \not\subseteq \gp$ for any $\gp \in \assh{R}{M}$,
and so $\assh{R}{M} = \pp \cap {\rm Supp}_{\,R}\,M$.
Thus we get the required equality.

\vspace{1em}

The proof of the next result was suggested by S. Goto.
\begin{thm}\label{3.8} \hspace{0.5ex}
If $0 \rightarrow L \rightarrow M \rightarrow N \rightarrow 0$
is an exact sequence of $R$-modules, then 
$\jmult{d}{I}{M} = \jmult{d}{I}{L} + \jmult{d}{I}{N}\,$.
\end{thm}
{\it Proof.}\hspace{0.5ex}
We use induction on $d$.
If $d = 1$, then $\length{R_\gp}{M_\gp} =
\length{R_\gp}{L_\gp} + \length{R_\gp}{N_\gp}$
for any $\gp \in {\rm Spec}\,R$ with $\kdim{}{R/\gp} = 1$,
and hence we get $\jmult{1}{I}{M} = \jmult{1}{I}{L}
+ \jmult{1}{I}{N}$ by \ref{3.7}.
So, let us consider the case where $d \geq 2$.
We set $\ic{L} = L/(0 :_L I^\infty)$,
$\ic{M} = M/(0 :_M I^\infty)$ and $\ic{N} = N/(0 :_N I^\infty)$.
Let $a$ be a sufficiently general element of $I$.
By \ref{3.4} it is enough to show
\begin{equation}\label{eq:1}
\jmult{d-1}{I}{\ic{M}/{a\ic{M}}} =
\jmult{d-1}{I}{\ic{L}/{a\ic{L}}} +
\jmult{d-1}{I}{\ic{N}/{a\ic{N}}}\,.
\end{equation}
Now we consider the exact sequences
\begin{eqnarray*}
 & & 0 \ra E \ra \ic{M} \ra \ic{N} \ra 0 \hspace{2ex}\mbox{and}  \\
 & & 0 \ra \ic{L} \ra E \ra F \ra 0\,,   
\end{eqnarray*}
where $E = (L :_M I^\infty)/(0 :_M I^\infty)$ and
$F = (L :_M I^\infty)/((0 :_M I^\infty) + L)$.
As $a$ can be chosen so that it is a non-zero-divisor on
$\ic{M}$ and $\ic{N}$, we get the exact seuences
\begin{eqnarray*}
 & & 0 \ra E/{aE} \ra \ic{M}/{a\ic{M}} \ra \ic{N}/{a\ic{N}}
    \ra 0  \hspace{2ex}\mbox{and}     \\
 & & 0 \ra {0 :_F a} \ra \ic{L}/{a\ic{L}} \ra E/{aE}
    \ra F/{aF} \ra 0 \,.   
\end{eqnarray*}
Then, the inductive hypothesis implies
\begin{eqnarray*}
\jmult{d-1}{I}{\ic{M}/{a\ic{M}}} & = &
  \jmult{d-1}{I}{E/{aE}} + \jmult{d-1}{I}{\ic{N}/{a\ic{N}}}  \\
\jmult{d-1}{I}{E/{aE}} & = &
  \jmult{d-1}{I}{X} + \jmult{d-1}{I}{F/{aF}} \hspace{2ex}
  \mbox{and}   \\
\jmult{d-1}{I}{\ic{L}/{a\ic{L}}} & = &
  \jmult{d-1}{I}{0 :_F a} + \jmult{d-1}{I}{X}\,,
\end{eqnarray*}
where $X$ denotes the kernel of $E/{aE} \rightarrow F/{aF}$.
Because $I^nF = 0$ for $n \gg 0$,
we have $\jmult{d-1}{I}{F/{aF}} = \jmult{d-1}{I}{0 :_F a} = 0$,
and so $\jmult{d-1}{I}{E/{aE}} = \jmult{d-1}{I}{X} =
\jmult{d-1}{I}{\ic{L}/{a\ic{L}}}$.
Therefore we get the required equality (\ref{eq:1})
and the proof is complete.
\begin{thm}\label{3.9} \hspace{0.5ex}
${\displaystyle
\jmult{d}{I}{M} = \sum_{\gp \in \assh{R}{M}}
\length{R_\gp}{M_\gp}\cdot \jmult{d}{I}{R/\gp}\,.}$
\end{thm}
{\it Proof.}\hspace{0.5ex}
Similarly as the proof of the additive formula for
the usual multiplicity, by \ref{3.8} we can prove
the required equality using induction on
$\sigma(M) =
\sum_{\gp \in \assh{R}{M}}\,\length{R_\gp}{M_\gp}$.
\begin{thm}\label{3.10}
Let $\ell(I) = \kdim{}{R} = d$ and $1 \leq i < d$.
Choosing general elements
$a_1, \dots\,, a_i$ of $I$,
we set $\ga = (a_1, \dots\,, a_i)$ and
$\pp = \{\,\gp \in \Min{R}{R / \ga} \mid
\mbox{$\height{R}{\gp} = i$ and $\kdim{}{R / \gp} = d - i$}\,\}$.
Then we have the following asertions.
\begin{itemize}
\item[{\rm (i)}]
$\pp \neq \phi$
\item[{\rm (ii)}]
$R_\gp$ is Cohen-Macaulay for any $\gp \in \pp$.
\item[{\rm (iii)}]
${\displaystyle
\jm{I} = \sum_{\gp \in \pp}\,
\jm{\ga R_\gp}\cdot\jm{I + \gp\,/\,\gp}\,.
}$
\end{itemize}
\end{thm}
{\it Proof.}\hspace{0.5ex}
Let $S = R / ((a_1, \dots\,, a_i) :_R I^\infty)$.
Then by \ref{3.5} and \ref{3.9} we have
$\kdim{R}{S} \leq d - i$ and
\[
\jm{I} = \sum_{\gp \in \assh{R}{S}}\,
\length{R_\gp}{S_\gp}\cdot\jmult{d-i}{I}{R / \gp}\,.
\]
Because $\jm{I} \neq 0$ by \ref{3.1},
there exists $\gp \in \assh{R}{S}$ such that
$\kdim{}{R / \gp} = d - i$.
Then $\height{R}{\gp} \leq i$,
and so $I \not\subseteq \gp$ by \cite[7.2]{aht}.
Therefore $S_\gp = R_\gp / \ga R_\gp$.
Here we notice that we can choose
$a_1, \dots\,, a_i$ so that
$a_j$ is regular on $R / ((a_1, \dots\,, a_{j-1}) :_R I^\infty)$
for any $1 \leq j \leq i$.
Then $a_1, \dots\,, a_i$ is a regular sequence on $R_\gp$,
and hence $R_\gp$ is a Cohen-Macaulay ring with
$\kdim{}{R_\gp} = i$.
Thus we have seen (i) together with (ii).
Furthermore we get (iii) since
$\pp = \{\,\gp \in \assh{R}{S} \mid
\kdim{}{R / \gp} = d - i\,\}$ and
$\jm{\ga R_\gp} = \mult{\ga R_\gp}{R_\gp} =
\length{R_\gp}{S_\gp}$ for any $\gp \in \pp$.
\section{Examples}
In this section we give some examples in order to
illustrate how to use \ref{3.6}.
\begin{ex}\label{4.1}
Let $S$ be a $3$-dimensional Cohen-Macaulay local ring
and $x, y, z$ be a system of parameters for $S$.
We set $R = S/{(x^2-yz)S}$ and $I = (x, y)R$.
Then we have
\[
\jm{I} = \length{S}{S/{(x, y, z)S}}\,.
\]
\end{ex}
{\it Proof.}\hspace{0.5ex}
If $\gp \in \ass{R}{R/I}$, then $\height{R}{\gp} = 1$
as $R/I = S/{(x, y)S}$ is a $1$-dimensional Cohen-Macaulay
ring. So the image of $z$ in $R$ becomes unit in $R_\gp$,
which means $IR_\gp = xR_\gp$.
Hence $I$ is generically a complete intersection.
On the other hand, if $f, g \in S$ satisfy
$xf - yg = (x^2 - yz)h$ for some $h \in S$,
we get $x(f - xh) = y(g - zh)$.
Hence there exists $\ell \in S$ such that
$f - xh = y\ell$ and $g - zh = x\ell$,
since $x, y$ form an $S$-regular sequence.
This means that we have an exact sequence
\begin{equation}
\label{4.2}
\begin{array}{ccccccc}
R^2 & \ra & R^2 & \ra & I & \ra & 0\,. \\
 & \left(\begin{array}{cc}
   y & x \\
   x & z
   \end{array}\right) & & \left(\begin{array}{cc}
                          x & -y
                          \end{array}\right) & & & 
\end{array}
\end{equation}
Now we choose a general element $\xi$ of $I$
so that $I = \xi R +yR$.
Then $\xi R :_R y = \xi R :_R I \subseteq \xi R :_R I^\infty$.
Let $\gq \in \ass{R}{R/(\xi R :_R y)}
\subseteq \ass{R}{R / {\xi R}}$.
Then $\height{R}{\gq} = 1$.
If $I \subseteq \gq$,
we have $IR_\gq = \xi R_\gq$
as $I$ is generically a complete intersection,
and so $(\xi R :_R I^\infty)R_\gq = R_\gq =
(\xi R :_R y)R_\gq$.
If $I \not\subseteq \gq$,
we have $(\xi R :_R I^\infty)R_\gq = \xi R_\gq =
(\xi R :_R y)R_\gq$
since $y \not\in \gq$.
Thus we get $\xi R :_R I^\infty = \xi R :_R y$,
and hence
\[
\jm{I} = \length{R}{R/{(\xi R :_R y) + yR}}
\]
by \ref{3.6} and \ref{3.7}.
Here we notice that $\xi$ can be written as the
image of $x + uy$ for some unit $u$ of $S$.
It is enough to show
$\xi R :_R y = (x + uy, z + ux)R$.
One can easily check that the right hand side is
contained in the left hand side.
Let us take any $\rho \in \xi R :_R y$
and write $y\rho = (x + uy)\eta$ for some $\eta \in R$.
Then $x\eta - y(\rho - u\eta) = 0$,
and so from the exact sequence (\ref{4.2}) we get
\[
\left(\begin{array}{c}
\eta \\
\rho - u\eta
\end{array}\right)
=
\left(\begin{array}{cc}
y & x \\
x & z
\end{array}\right)\,
\left(\begin{array}{c}
\sigma \\
\tau
\end{array}\right)
\]
for some $\sigma, \tau \in R$.
This means
$\rho = \sigma (x + uy) + \tau (z + ux)$.
Thus we see $\xi R :_R y = (x + uy, z + ux)R$
and the proof is complete.
\begin{ex}\label{4.3}
Let $(R, \gm)$ be a $3$-dimensional Cohen-Macaulay local ring
and $I$ an ideal of height two which is generated by
the maximal minors of the matrix
\[
\varphi = \left(\begin{array}{lll}
x_{11} & x_{12} & x_{13} \\
x_{21} & x_{22}& x_{23}
\end{array}\right)
\]
with entries in $R$.
We set $S = R/I$ and asume that each one of
$x_{11}, x_{12}$ and $x_{13}$ is a system of parameters for $S$.
\begin{itemize}
\item[{\rm (i)}]
For some units $u_1, u_2$ and $u_3$ of $R$, we have 
\[
\jm{I} = \length{S}{
S/{(x_{11}u_1+x_{12}u_2+x_{13}u_3, x_{21}u_1+x_{22}u_2+x_{23}u_3)S}}\,.
\]

\item[{\rm (ii)}]
Suppose that $R$ is a Gorenstein ring.
Then $(x_{1j}, x_{2j})S$ is a canonical ideal of $S$
for $j = 1, 2, 3$.
On the other hand, any canonical ideal of $S$ is
written in the form
$(x_{11}a_1+x_{12}a_2+x_{13}a_3, x_{21}a_1+x_{22}a_2+x_{23}a_3)S$
for some elements $a_1, a_2$ and $a_3$ of $R$.
\item[{\rm (iii)}]
Suppose that $S$ is an integral domain
whose normalization is a DVR.
Then we have
\[
\jm{I} = \min\{\,
\length{S}{S/{(x_{1j}, x_{2j})S}}\, \}_{j = 1, 2, 3}\,.
\]
Moreover, if $R$ is a Gorenstein ring, we get
\[
\jm{I} = \min\{\,
\length{S}{S/\omega} \mid
\mbox{$\omega$ is a canonical ideal of $S$}\, \}\,.
\]
\end{itemize}
\end{ex}
{\it Proof.}\hspace{1ex}
Let $f_j$ be the determinant of the $2 \times 2$ matrix
derived from $\varphi$ deleting the $j$-th column.
It is well known that
\begin{equation}\label{4.4}
\begin{array}{ccccccccc}
0 & \ra & R^2 & \ra & R^3 & \ra & I \ra & 0 \\
  &  &  & \trmatrix{\,\varphi} &  &
  (f_1\hspace{1ex}{-f_2}\hspace{1ex}f_3) & & &
\end{array}
\end{equation}
is an exact sequence.
In particular, we have that $R/I$ is a Cohen-Macaulay ring and
\begin{equation}\label{4.5}
x_{i1}f_1 - x_{i2}f_2 + x_{i3}f_3 = 0
\hspace{2ex}\mbox{for $i =1, 2$}\,.
\end{equation}
If $\gp \in \assh{R}{R/I}$,
then each one of $x_{11}, x_{12}$ and $x_{13}$
is a unit in $R_\gp$, and so the equality (\ref{4.5})
for $i = 1$ implies that
$IR_\gp$ is generated by any two of $f_1, f_2$ and $f_3$.
Hence $I$ is generically a complete intersection.
On the other hand, for any $j = 1, 2, 3$,
we have that $x_{1j}, x_{2j}, f_j$ is
a system of parameters for $R$ as
$(x_{1j}, x_{2j}, f_j) \supseteq (x_{1j}) + I$,
and so $x_{1j}, x_{2j}$ is an $R$-regular sequence.
Moreover, we see that any two of
$f_1, f_2$ and $f_3$ form an $R$-regular sequence.
In fact, for example,
if there exists a height one prime ideal $\gq$
of $R$ containing $(f_1, f_2)$,
then $f_3 \not\in \gq$ as $\height{R}{I} = 2$,
and so from (\ref{4.5}) we get
$(x_{13}, x_{23}) \subseteq \gq$,
which is impossible as $x_{13}, x_{23}$
is an $R$-regular sequence.

(i)\hspace{1ex}
Let $g_1$ and $g_2$ be sufficiently general elements of $I$.
Then
\[
\jm{I} = \length{S}{S/{((g_1, g_2) :_R I^\infty})S}
\]
by \ref{3.6}.
Let $\gp \in \ass{R}{R/((g_1, g_2) :_R I)} \subseteq
\ass{R}{R/(g_1, g_2)}$.
Then $\height{R}{\gp} = 2$ since $g_1, g_2$
is an $R$-regular sequence.
Therefore, if $I \subseteq \gp$,
we have $IR_\gp = (g_1, g_2)R_\gp$, and so
$((g_1, g_2) :_R I^\infty)R_\gp =
((g_1, g_2) :_R I)R_\gp = R_\gp$.
Even if $I \not\subseteq \gp$, we have
$((g_1, g_2) :_R I^\infty)R_\gp =
((g_1, g_2) :_R I)R_\gp = (g_1, g_2)R_\gp$.
Thus we get $(g_1, g_2) :_R I^\infty
= (g_1, g_2) :_R I$.
We write
\[
g_1 = v_1f_1 + v_2f_2 + v_3f_3 \hspace{2ex}
\mbox{and} \hspace{2ex}
g_2 = w_1f_1 + w_2f_2 + w_3f_3\,,
\]
where $v_j$ and $w_j$ are units of $R$.
We may assume that all maximal minors of the matrix
\[
\psi = \left(\begin{array}{lll}
v_1 & v_2 & v_3 \\
w_1 & w_2 & w_3
\end{array}\right)
\]
are units of $R$.
We set
\[
h_1 = v_2g_2 - w_2g_1 \hspace{2ex}
\mbox{and} \hspace{2ex}
h_2 = v_1g_2 - w_1g_1\,.
\]
Then $(g_1, g_2) = (h_1, h_2)$.
Moreover, we have
\[
h_1 = -u_3f_1 + u_1f_3 \hspace{2ex}
\mbox{and} \hspace{2ex}
h_2 = u_3f_2 + u_2f_3\,,
\]
where $u_j$ denotes the determinant of the matrix
derived from $\psi$ deleting the $j$-th column.
Hence we get $(h_1, h_2, f_3) = I$, and so
$(g_1, g_2) :_R I = (h_1, h_2) :_R f_3$.

Let us take any $\rho \in (h_1, h_2) :_R f_3$.
Then there exist $\xi, \eta \in R$ such that
\begin{eqnarray*}
\rho f_3 & = & \xi h_1 + \eta h_2 \\
         & = & \xi (-u_3 f_1 + u_1f_3) + \eta (u_3f_2 + u_2f_3)\,,
\end{eqnarray*}
which induces
\[
\xi u_3 f_1 + \eta u_3(-f_2) + (\rho - \xi u_1 - \eta u_2) f_3 = 0\,.
\]
Therefore, using the exact sequence (\ref{4.4}),
we get
\[
\left(\begin{array}{c}
\xi u_3 \\
\eta u_3 \\
\rho - \xi u_1 - \eta u_2
\end{array}\right)
=
\sigma
\left(\begin{array}{c}
x_{11} \\
x_{12} \\
x_{13}
\end{array}\right)
+ \tau
\left(\begin{array}{c}
x_{21} \\
x_{22} \\
x_{23}
\end{array}\right)
\]
for some $\sigma, \tau \in R$.
This implies $u_3\rho = \sigma y_1 + \tau y_2$,
where $y_i = x_{i1}u_1 + x_{i2}u_2 + x_{i3}u_3$
for $i = 1, 2$.
Thus we see $(h_1, h_2) :_R f_3 \subseteq (y_1,\, y_2)$.
On the other hand, the converse inclusion holds since
$y_3f_3 = x_{i1}h_1 + x_{i2}h_2$
for $i = 1, 2$,
and so the proof of (i) is complete.

(ii)\hspace{1ex}
Let us prove that $(x_{11}, x_{21})S$
is a canonical ideal of $S$.
By the same argument one can prove that
$(x_{12}, x_{22})S$ and $(x_{13}, x_{23})S$
are also canonical ideals of $S$.

Let $\gp \in \ass{R}{R/(f_2, f_3)}$.
Then $\height{R}{\gp} = 2$ as
$f_2, f_3$ is an $R$-regular sequence.
Hence we have $IR_\gp = (f_2, f_3)R_\gp$
when $I \subseteq \gp$.
On the other hand, if $I \not\subseteq \gp$,
$f_1$ is a unit of $R_\gp$, and so
$(x_{11}, x_{21})R_\gp = (f_2, f_3)R_\gp$ by (\ref{4.5}).
Thus we get $(x_{11}, x_{21}) \cap I = (f_2, f_3)$.
Then
\begin{eqnarray*}
(x_{11}, x_{21})S & \cong &
       (x_{11}, x_{21})\, /\, {(x_{11}, x_{21}) \cap I} \\
 & = & (x_{11}, x_{21})\, /\, (f_2, f_3)  \\
 & = & ((f_2, f_3) :_R f_1)\, /\, (f_2, f_3) \\
 & \cong & {\rm Hom}_{\,R}(S, R / (f_2, f_3))\,.
\end{eqnarray*}
Because $S$ is a homomorphic image of the Gorenstein local
ring $R/(f_2, f_3)$, it follows that
${\rm Hom}_{\,R}(S, R / (f_2, f_3))$ is the canonical module of $S$,
and so $(x_{11}, x_{21})S$ is a canonical ideal of $S$.

Now let $\omega$ be any canonical ideal of $S$.
Let $Q$ be the total quotient ring of $S$.
Then there exists $\alpha \in Q$ such that
$\omega = \alpha \cdot (x_{11}, x_{21})S$.
Because $S :_Q (x_{11}, x_{21})S =
{\ic{x_{11}}}^{\,-1}(x_{11}S :_S x_{21})$,
where $\ic{\,\ast\,}$ denotes the reduction mod $I$,
we write $\alpha = {\ic{x_{11}}}^{\,-1}\cdot\ic{y}$
with $y \in R$ such that
$\ic{y} \in x_{11}S :_S x_{21}$.
Then $x_{21}y \in (x_{11}) + I =
(x_{11}, x_{12}x_{21}, x_{13}x_{21}, f_1)$,
and so there exist $z, z' \in R$ such that
$x_{21}(y - x_{12}z - x_{13}z') \in (x_{11}, f_1)$.
This means $y - x_{12}z - x_{13}z' \in (x_{11}, f_1)$
since $x_{11}, x_{21}, f_1$ is an $R$-regular sequence.
Hence $y \in (x_{11}, x_{12}, x_{13}, f_1) =
(x_{11}, x_{12}, x_{13})$. Therefore
$y = x_{11}a_1+ x_{12}a_2 + x_{13}a_3$ for some
$a_1, a_2, a_3 \in R$. Then we have
\begin{eqnarray*}
\omega & = & 
       {\ic{x_{11}}}^{\,-1}\cdot\sum_{j=1}^3
         \ic{x_{1j}}\,\ic{a_j} \cdot (x_{11},\, x_{21})S  \\
 & = & (\,\sum_{j=1}^3 \ic{x_{1j}}\,\ic{a_j}\,, \,
        \sum_{j=1}^3
        {\ic{x_{11}}}^{\,-1}\ic{x_{1j}}\,\ic{x_{21}}\,\ic{a_j}\,)\,.
\end{eqnarray*}
Because ${\ic{x_{11}}}^{\,-1}\ic{x_{1j}}\,\ic{x_{21}} =
\ic{x_{2j}}$ for $j = 1, 2, 3$, we get
\[
\omega = (\,\sum_{j=1}^3 x_{1j}a_j\,, \,\sum_{j=1}^3 x_{2j}a_j\,)S
\]
and the proof of (ii) is complete.

In order to prove (iii), we need the following lemma.
\begin{lem}\label{4.6}
Let $(S, \gn)$ be a $1$-dimensional Cohen-Macaulay local ring
and $T$ be an extension ring of $S$ such that
$T / S$ has finite length as an $S$-module.
Let $\ga$ be an $\gn$-primary ideal of $S$ and $v$
be a non-zero-divisor of $T$
such that $\ga v \subseteq S$.
Then $\ga v$ is an $\gn$-primary ideal of $S$ and
\[
\length{S}{S / {\ga v}} \geq \length{S}{S / \ga}\,,
\]
where the equality holds if $v$ is a unit of $T$.
\end{lem}
{\it Proof.}\hspace{1ex}
Let $\gc = S :_S T$,
which is an $\gn$-primary ideal of $S$.
We take an $S$-regular element $u \in \ga$
and set $\gb = u\gc$,
which is an $\gn$-primary ideal of $S$ contained in $\ga$.
Let us consider the exact sequence
\[
0 \ra {\ga v} / {\gb v} \ra S / {\gb v} \ra S / {\ga v} \ra 0\,.
\]
As $v$ is a non-zero-divisor of $T$,
we have ${\ga v} / {\gb v} \cong \ga / \gb$, and so
\[
\length{S}{S / {\ga v}} =
\length{S}{S / {\gb v}} - \length{S}{\ga / \gb}\,.
\]
Now we notice $\gb v = (u\gc)v = u(\gc v)
\subseteq u\gc = \gb$.
Moreover, if $v$ is a unit of $T$,
$\gb = (u\gc)v^{-1}v = u(\gc v^{-1})v \subseteq
u\gc v = \gb v$.
Thus we get the required assertion.

\vspace{1em}
\noindent
{\it Proof of (iii) of \ref{4.3}.}\hspace{1ex}
Let $T$ be the normalization of $S$.
Let $tT$ be the maximal ideal of $T$.
We write $\ic{x_{ij}} = t^{k_{ij}}\alpha_{ij}$,
where $0 \leq k_{ij} \in \zz$ and
$\alpha_{ij}$ is a unit of $T$.
Replacing the columns of $\varphi$,
we may assume $k_{11} \leq k_{12}$ and $k_{11} \leq k_{13}$.

Because $\ic{x_{11}}\,\ic{x_{22}} - \ic{x_{12}}\,\ic{x_{21}} = 0$ 
and
$\ic{x_{11}}\,\ic{x_{23}} - \ic{x_{13}}\,\ic{x_{21}} = 0$,
we have
\begin{equation}\label{4.7}
t^{k_{11}+k_{22}}\alpha_{11}\alpha_{22} =
t^{k_{12}+k_{21}}\alpha_{12}\alpha_{21}
\hspace{2ex}\mbox{and}\hspace{2ex}
t^{k_{11}+k_{23}}\alpha_{11}\alpha_{23} =
t^{k_{13}+k_{21}}\alpha_{13}\alpha_{21}\,.
\end{equation}
Now we put
$\beta_j = t^{k_{1j} - k_{11}} {\alpha_{11}}^{-1} \alpha_{1j}
\in T$ for $j = 2, 3$.
Then by (\ref{4.7}) we have
$\ic{x_{ij}} = \ic{x_{i1}}\cdot \beta_j$
for any $i = 1, 2$ and $j = 2, 3$,
and so $(x_{1j}, x_{2j})S =
\beta_j \cdot (x_{11}, x_{21})S$ for $j = 2, 3$.
Therefore by \ref{4.6} we get
$\length{S}{S / (x_{11}, x_{21})S} \leq
\length{S}{S / (x_{1j}, x_{2j})S}$ for $j = 2, 3$.

Next we prove
$\jm{I} = \length{S}{S / (x_{11}, x_{21})S}$.
Let $u_1, u_2$ and $u_3$ be units of $S$ stated in (i).
Then $\jm{I} = \length{S}{S / (y_1, y_2)S}$,
where $y_i = x_{i1}u_1 + x_{i2}u_2 + x_{i3}u_3$.
Here we put $\gamma = \ic{u_1} +
t^{k_{12} - k_{11}}{\alpha_{11}}^{-1}\alpha_{12}\ic{u_2} +
t^{k_{13} - k_{11}}{\alpha_{11}}^{-1}\alpha_{13}\ic{u_3}$.
Then $\ic{y_1} = \ic{x_{11}}\gamma$
and $\ic{y_2} = \ic{x_{12}}\gamma$
by (\ref{4.7}).
Let us notice that the units $v_j$ and $w_j$
stated in the proof of (i) can be chosen so that
$\gamma$ is a unit in $T$.
Hence we get $\length{S}{S / (y_1, y_2)S} =
\length{S}{S / (x_{11}, x_{21})S}$ by \ref{4.6}.

The last assertion of (iii) can be proved by
the similar argument using (ii),
and so the proof of \ref{4.3} is complete.
\begin{ex}\label{3c}
Let $R = K[[X, Y, Z]]$ and $T = K[[t]]$
be formal power series rings over an infinite field $K$.
Let $\varphi : R \ra T$ be the homomorphism of $K$-algebras
defined by $\varphi(X) = t^k$, $\varphi(Y) = t^\ell$,
$\varphi(Z) = t^m$, where $k, \ell, m$ are positive integers
such that ${\rm GCD}\{ k, \ell, m \} = 1$.
We put $\gp = {\rm Ker}(\varphi)$.
Then $\gp$ is generated by the maximal minors of the matrix
of the form
\[
\left(\begin{array}{lll}
X^\alpha & Y^{\beta'} & Z^{\gamma'} \\
Y^\beta  & Z^\gamma   & X^{\alpha'}
\end{array}\right)\,,
\]
where $\alpha, \beta, \gamma, \alpha', \beta', \gamma'$
are positive integers {\rm (cf. \cite{h})}.
By replacing the variables X, Y and Z suitably,
we may assume
$k\alpha = \min\{ k\alpha, \ell\beta, m\gamma,
k\alpha', \ell\beta', m\gamma' \}$.
Then we have
\[
\jm{\gp} = \alpha\beta(\gamma + \gamma')\,.
\]
\end{ex}
{\it Proof.}\hspace{0.5ex}
We put $S = R / \gp = K[[t^k, t^\ell, t^m]]$.
By (iii) of \ref{4.3}, 
$\jm{\gp}$ is the minimum of the lengths of
$S / (t^{k\alpha}, t^{\ell\beta})S$,
$S / (t^{\ell\beta'}, t^{m\gamma})S$ and
$S / (t^{m\gamma'}, t^{k\alpha'})S$.
Here we notice that $\ell(\beta + \beta') = k\alpha + m\gamma$
as $Y^{\beta + \beta'} - X^\alpha Z^\gamma \in \gp$.
This means $m\gamma = \ell\beta' - k\alpha + \ell\beta$,
and so $(t^{\ell\beta'}, t^{m\gamma})S =
t^\delta \cdot (t^{k\alpha}, t^{\ell\beta})S$,
where $\delta = \ell\beta' - k\alpha \geq 0$.
Hence we have $\length{S}{S / {(t^{k\alpha}, t^{\ell\beta})S}}
\leq \length{S}{S / {(t^{\ell\beta'}, t^{m\gamma})S}}$
by \ref{4.6}.
Similary we get
$\length{S}{S / {(t^{k\alpha}, t^{\ell\beta})S}} \leq
\length{S}{S / {(t^{m\gamma'}, t^{k\alpha'})S}}$.
Therefore $\jm{\gp} = \length{S}{S / {(t^{k\alpha}, t^{\ell\beta})S}}$.
Because $S / {(t^{k\alpha}, t^{\ell\beta})S} \cong
R / {\gp + (X^\alpha, Y^\beta)R} =
R / {(X^\alpha, Y^\beta, Z^{\gamma + \gamma'})R}$,
we get the required equality.
\vspace{1em}

\begin{ex}\label{3d}
Let $S, J$ and $I$ be as in \ref{3a}.
Moreover we assume that ${\rm Proj}(S)$ is Cohen-Macaulay
and $JS_P$ is weakly $(d - 3)$-residually $S_2$ for any
$P \in {\rm V}(J) \setminus \{\,\gn\,\}$.
Then, if $\gb \subseteq J$ is an ideal generated by
$d - 1$ homogeneous elements of degree $r$ such that
$\height{S}{(\gb :_S J)} \geq d - 1$,
\[
\jm{I} = r \cdot \mult{}{S / (\gb :_S J)}\,.
\]
\end{ex}
{\it Proof.}\hspace{0.5ex}
By \ref{3a} we have $\jm{I} = r \cdot \mult{}{S / {\ga :_S J}}$,
where $\ga$ is an ideal generated by $d - 1$ {\it general} homogeneous
elements of $J$ of degree $r$.
On the other hand, by \cite[2.1 and 3.1]{u} we have
$\mult{}{S / (\ga :_S J)} = \mult{}{S / (\gb :_S J)}$.
Therefore we get the required equality.
\begin{ex}\label{3e}
Let $S = \oplus_{n \geq 0}\,S_n$ be a $d$-dimensional
Gorenstein standard algebra such that $S_0$ is an
infinite field.
Let $I$ be a homogeneous ideal of 
generically a complete intersection such that
$S / I$ is a $2$-dimensional Cohen-Macaulay ring.
Let $J$ be a reduction of $I$ generated by homogeneous
elements $g_1, \dots\,, g_{d-1}, g_d$ of the same degree.
We set $\gn = S_{+}$ and $\gb = (g_1, \dots\,, g_{d-1})$.
Then, if $(\gb :_S J) + J$ is $\gn$-primary,
we have
\[
\jm{IS_\gn} = \length{S}{S / {(\gb :_S J) + J}}\,.
\]
\end{ex}
{\it Proof.}\hspace{0.5ex}
First we show that $J$ satisfies $G_d$ and
$J$ is $(d - 2)$-residually $S_2$.
We take any $P \in {\rm V}(J) = {\rm V}(I)$
with $\height{S}{P} \leq d - 1$.
Then $\gb :_S J \not\subseteq P$, and so
$\gb S_P = JS_P$.
In particular, $\mu_{S_P}(JS_P) \leq d - 1$.
If $\height{S}{P} = d - 2$,
$P \in \Min{S}{S / I}$,
and so $IS_P$ is a complete intersection,
which means $IS_P = JS_P$ as $JS_P$ is a
reduction of $IS_P$.
Hence $\mu_{S_P}(JS_P) = d - 2$ if $\height{S}{P} = d - 2$.
Hence $J$ satisfies $G_d$.
In order to prove that $J$ is $(d - 2)$-residually $S_2$,
it is enough to show that if $\gc \subseteq J$ is
$(d - 2)$-generated ideal with $\height{S}{(\gc :_S J)} \geq d - 2$,
then $S / (\gc :_S J)$ is Cohen-Macaulay.
Because $\height{S}{(\gc :_S J)} \geq d - 2$ and
$\height{S}{J} = d - 2$, we have $\height{S}{\gc} = d - 2$,
and so $\gc$ is generated by a regular sequence.
We take any $Q \in \ass{S}{S / (\gc :_S I)}$.
Then $\height{S}{Q} = d - 2$ as $S / (\gc :_S I)$ is a
$2$-dimensional Cohen-Macaulay ring (cf. \cite[1.3]{ps}).
If $I \subseteq Q$, then $Q \in \Min{S}{S / I}$,
and so $IS_Q = JS_Q$,
which means $(\gc :_S I)S_Q = (\gc :_S J)S_Q$.
On the other hand, if $I \not\subseteq Q$, we have
$J \not\subseteq Q$, and so
$(\gc :_S I)S_Q = \gc S_Q = (\gc :_S J)S_Q$.
Consequently, $\gc :_S I = \gc :_S J$,
and hence $S / (\gc :_S J)$ is Cohen-Macaulay.

Therefore, by \ref{3b} and \ref{3.11}
we have
\[
\jm{JS_\gn} =
\length{S}{S / ((f_1, \dots\,, f_{d-1}) :_S J) + J}
\]
for {\it general} homogeneous elements
$f_1, \dots\,, f_{d-1} \in J$ of degree $r$.
We put $\ga = (f_1, \dots\,, f_{d-1})$.
Then, as $\ga :_S J$ is a geometric $(d - 1)$-residual
intersection of $J$, we have
$(\ga :_S J) \cap J = \ga$ by \cite[3.4]{ceu}.
Thus we get
\begin{equation}\label{4.10}
0 \ra S / \ga \ra S / (\ga :_S J) \oplus S / J
\ra S / {(\ga :_S J) + J} \ra 0
\end{equation}
of graded $S$-modules.
Similary, as $\gb :_S J$ is also a geometric
$(d - 1)$-residual intersection of $J$,
we have the exact sequence
\begin{equation}\label{4.11}
0 \ra S / \gb \ra S / (\gb :_S J) \oplus S / J
\ra S / {(\gb :_S J) + J} \ra 0
\end{equation}
of graded $S$-modules.
Here we notice that by \cite[2.2]{ceu}
the Hilbert function of $S / \ga$ (resp. $S / (\ga :_S J)$)
coincides with that of $S / \gb$ (resp. $S / (\gb :_S J)$).
Therefore by (\ref{4.10}) and (\ref{4.11}) we see that
the Hilbert functions of $S / {(\ga :_S J) + J}$ and
$S / {(\gb :_S J) + J}$ are same,
and hence the lengths of those algebras are same.
Therefore the proof is complete as
$\jm{IS_\gn} = \jm{JS_\gn}$ by \cite[2.10]{fm2}
\begin{ex}\label{3f}
Let $k$ be an infinite field.

\noindent
\begin{itemize}
\item[{\rm (1)}]
Let $S = k[X_1, X_2, X_3, X_4, X_5]$ be a polynomial ring and
$I$ be the ideal generated by the maximal minors of the matrix
\[
\left(\begin{array}{llll}
X_1 & X_2 & X_3 & X_4 \\
X_2 & X_3 & X_4 & X_5
\end{array}\right)\,.
\]
We set $\gn = S_+$. Then $\jm{IS_\gn} = 4$.
\item[{\rm (2)}]
Let $S = k[X_1, X_2, X_3, X_4, X_5, X_6, X_7]$
be a polynomial ring and $I$ be the ideal generated
by the maximal minors of the matrix
\[
\left(\begin{array}{lllll}
X_1 & X_2 & X_3 & X_5 & X_6 \\
X_2 & X_3 & X_4 & X_6 & X_7
\end{array}\right)\,.
\]
We set $\gn = S_+$.
Then $\jm{IS_\gn} = 10$.
\end{itemize}
\end{ex}
{\it Proof.}\hspace{1ex}
(1)\hspace{0.5ex}
It is well known that $I$ is a perfect ideal
of height $3$ which is generically a complete intersection.
Let $\Delta_{ij}$ be the determinant of the matrix
\[
\left(\begin{array}{ll}
X_i & X_j \\
X_{i+1} & X_{j+1}
\end{array}\right)
\]
for $1 \leq i < j \leq 4$.
We set $\ga = (\Delta_{12} - \Delta_{34},
\Delta_{13}, \Delta_{14}, \Delta_{24})$ and
$J = \ga + \Delta_{23}R$.
Then $J$ is a reduction of $I$.
Moreover, we have
$\ga :_R \Delta_{23} = (X_2, X_4, X_5 - X_1, {X_1}^2, X_1X_5)$,
and so
$(\ga : \Delta_{23}) + \Delta_{23}R =
(X_2, X_4, X_5 - X_1, {X_1}^2, X_1X_5, {X_3}^2)$,
which is an $\gn$-primary ideal.
Hence $\jm{IS_\gn} =
\length{S}{S / (\ga :_S \Delta_{23}) + \Delta_{23}S} = 4$
by \ref{3e}.

(2) This assertion can be proved by \ref{3d}.

\end{document}